\documentclass[12pt,draft]{article}

\usepackage{amsfonts,amsmath,amssymb,theorem}

\usepackage[english]{babel}

\DeclareMathOperator{\Id}{Id}
\DeclareMathOperator{\ap}{ap}
\DeclareMathOperator{\adj}{adj}
\DeclareMathOperator{\cp}{cap}
\DeclareMathOperator{\Lip}{Lip}
\DeclareMathOperator{\ACL}{ACL}

\DeclareMathOperator{\loc}{loc}
\DeclareMathOperator{\const}{const}

\DeclareMathOperator{\dist}{dist}
\DeclareMathOperator{\diam}{diam} 
\DeclareMathOperator{\esssup}{esssup}
\DeclareMathOperator{\ess}{ess}
 
\DeclareMathOperator{\BVL}{BVL} 

\begin{document}

\centerline{\bf Sobolev Homeomorphisms and Composition Operators
\footnote{This work was partially supported by Israel Scientific Foundation (Grant 1033/07)}}
\vskip 0.5cm

\centerline{\bf V.~Gol'dshtein and A.~Ukhlov}
\vskip 0.5cm

\centerline{ABSTRACT} 
\bigskip
{\small  We study invertibility of bounded composition operators of Sobolev spaces. The problem is closely connected with the theory of mappings of finite distortion.
If a homeomorphism $\varphi$ of Euclidean domains $D$ and $D'$ generates by the composition rule $\varphi^{\ast}f=f\circ\varphi$ a bounded composition operator of Sobolev spaces $\varphi^{\ast}: L^1_{\infty}(D')\to L^1_p(D)$, $p>n-1$, has finite distortion and Luzin $N$-property then its inverse $\varphi^{-1}$ generates the bounded composition operator from $L^1_{p'}(D)$, $p'=p/(p-n+1)$, into $L^1_{1}(D')$.  

 }
\vskip 0.5cm

\centerline{\bf Introduction}

\bigskip

Let $\varphi$ be a homeomorphism of Euclidean domains $D,D'\subset\mathbb R^n$. It is known [1] that $\varphi$ is a quasiconformal mapping if and only if the composition operator $\varphi^{\ast}$ is an isomorphism of Sobolev spaces $L^1_n(D')$ and $L^1_n(D)$. 
If $\varphi$ generates a bounded composition operator of Sobolev spaces $L^1_q(D')$ and $L^1_q(D)$, $q\ne n$, then the inverse homeomorphism $\varphi^{-1}$ is not necessary generates the bounded composition operator of same spaces. In the more general case homeomorphisms that generate composition operators 
from $L^1_p(D')$ to $L^1_q(D)$, $1\leq q\leq p\leq\infty$, are mappings with bounded $(p,q)$-distortion. These classes of mappings were introduced in [2] as a natural solution of the change of variable problem in Sobolev spaces. Inverse mappings to homeomorphisms with bounded $(p,q)$-distortion can be described in the same category of mappings with bounded mean distortion.  In [3] these classes of  mappings were studied in a relation with Sobolev type embedding theorems for non-regular domains.

We recall, that Sobolev space
$L^1_p(D)$, $1\leq p\leq\infty$, consists of locally summable, weakly differentiable functions $f:D\to\mathbb R$ with the finite seminorm:
$$
\|f\mid L^1_p(D)\|=\| |\nabla f|\mid L_p(D)\|,\quad \nabla f = \bigl(\frac{\partial f}{\partial x_1},...,\frac{\partial f}{\partial x_n}\bigr).
$$

As usually Lebesgue space $L_p(D)$, $1\leq p\leq\infty$, is the space of locally summable functions with the finite norm:
$$
\| f\mid L_p(D)\|=\biggr(\int\limits_D |f|^p~dx\biggl)^{\frac{1}{p}},\quad 1\leq p<\infty,
$$
and
$$
\|f\mid L_{\infty}(D)\|=
\ess\sup\limits_{x\in D}|f(x)|,\quad p=\infty.
$$

A mapping $\varphi:D\to\mathbb R^n$ belongs to $L^1_p(D)$, $1\leq p\leq\infty$, if
its coordinate functions $\varphi_j$ belong to $L^1_p(D)$, $j=1,\dots,n$.
In this case formal Jacobi matrix
$D\varphi(x)=\bigl(\frac{\partial \varphi_i}{\partial x_j}(x)\bigr)$, $i,j=1,\dots,n$,
and its determinant (Jacobian) $J(x,\varphi)=\det D\varphi(x)$ are well defined at
almost all points $x\in D$. The norm $|D\varphi(x)|$ of the matrix
$D\varphi(x)$ is the norm of the corresponding linear operator $D\varphi (x):\mathbb R^n \rightarrow \mathbb R^n$ defined by the matrix $D\varphi(x)$. We will use the same notation for this matrix and the corresponding linear operator.

Recall that a mapping $\varphi: D\to D'$ is called a {\it mapping
with bounded $(p,q)$-distortion} $1\leq q\leq p\leq\infty$, if $\varphi$ belongs to Sobolev
space $W^1_{1,\loc}(D)$ and the local $p$-distortion
$$
K_p(x)=\inf\{ k: |D\varphi|(x)\leq k |J(x,\varphi)|^{\frac{1}{p}}, \,\,x\in D\}
$$
belongs to Lebesgue space $L_{r}(D)$, where ${1}/{r}={1}/{q}-{1}/{p}$ (if $p=q$ then $r=\infty$).

Mappings with bounded $(p,q)$-distortion have
a finite distortion, i.~e.  $D\varphi(x)=0$ for almost all points $x$
that belongs to set $Z=\{x\in D:J(x,\varphi)=0\}$.

Necessity of studying of Sobolev mappings with integrable distortion arises in problems of the non-linear elasticity theory [4,\,5].
In these works J.~M.~Ball introduced classes of mappings, defined on bounded domains $D \in \mathbb R^n$:
$$
A^+_{p,q}(D)=\{\varphi\in W^1_p(D) : \adj D\varphi\in L_q(D),\quad J(x,\varphi)>0 \quad\text{a.~e. in}\quad D\},
$$
$p,q>n$, where $\adj D\varphi$ is the formal adjoint matrix to the Jacobi matrix $D\varphi$:
$$
\adj D\varphi(x)\cdot D\varphi(x) = \Id J(x,\varphi).
$$

The class of mappings with bounded $(p,q)$-distortion is a natural generalization of mappings with bounded distortion and represents a non-homeomorphic case of so-called $(p,q)$-quasi\-con\-for\-mal mappings [2,\,3,\,6,\,7]. Such classes of mappings have applications to the Sobolev type embedding problems [7--9].

The following assertion demonstrates a connection between Sobolev spaces and mappings with bounded $(p,q)$-distortion [2]. A homeomorphism $\varphi$ of Euclidean domains $D$ and $D'$ is a mapping with bounded $(p,q)$-distortion, $1\leq q\leq p<\infty$, if and only if $\varphi$ generates a bounded operator of Sobolev spaces
$$
\varphi^{\ast}: L^1_p(D')\to L^1_q(D)
$$
by the composition rule $\varphi^{\ast}f=f\circ\varphi$. We call $\varphi^{\ast}$ a composition operator of Sobolev spaces.

In the frameworks of the inverse operator problem in [6] was proved, that if a homeomorphism $\varphi : D\to D^{\prime}$ generates a bounded composition operator
$$
\varphi^{\ast}: L^1_{p}(D')\to L^1_q(D),\,\,\,n-1<q\leq p<+\infty,
$$
then the inverse mapping $\varphi^{-1}: D'\to D$ generates a bounded composition operator
$$
(\varphi^{-1})^{\ast}: L^1_{q'}(D)\to L^1_{p'}(D'),\,\,\,q'=q/(q-n+1),\,\,p'=p/(p-n+1).
$$

The main result of the article concerns to invertibility of a composition operator in the limit case $p=\infty$.

\vskip 0.3 cm

{\bf Theorem A.}
Let a homeomorphism $\varphi : D\to D^{\prime}$ has finite distortion, Luzin $N$-property (the image of a set measure zero is a set measure zero) and  generates a bounded composition operator
$$
\varphi^{\ast}: L^1_{\infty}(D')\to L^1_q(D),\,\,\,q>n-1.
$$
Then the inverse mapping $\varphi^{-1}: D'\to D$ generates a bounded composition operator
$$
(\varphi^{-1})^{\ast}: L^1_{q'}(D)\to L^1_1(D'),\,\,\,q'=q/(q-n+1).
$$
\vskip 0.3 cm

The invertibility problem for composition operators in Sobolev spaces is closely connected with a regularity problem for invertible Sobolev mappings.
The regularity problem for mappings which are inverse to Sobolev
homeomorphisms was studied by many authors. In article [10]
was proved that if a mapping $\varphi\in W^1_{n,\loc}(D)$ and
$J(x,\varphi)>0$ for almost all points $x\in D$, then $\varphi^{-1}$ belongs
to $W^1_{1,\loc}(D')$.

\medskip
 
{\em The assumption that $\varphi$ has finite distortion cannot be dropped out}. Indeed, consider the function $g(x)=x+u(x)$
on the real line, where $u$ is the standard Cantor function. Let
$f=g^{-1}$. Then the derivative $f'=0$ on the set of positive measure and $h^{-1}$  fails to be absolutely continuous. In this case we can prove only that the inverse
homeomorphism has a finite variation on almost all lines [11]. 
In work [11] was obtained the following result: {\em if a
homeomorphism $\varphi:D\to D'$ belongs to the Sobolev space
$L^1_p(D)$, $p>n-1$, then the inverse mapping $\varphi^{-1}: D'\to D$ has a finite variation on almost all lines (belongs to $\BVL(D')$)}.

\medskip

In work [12] the local regularity of plane homeomorphisms
that belong to Sobolev space $W^1_1(D)$ was studied. For the case of space $\mathbb R^n$, $n\geq 3$, recent work [13] contains the following result for domains in $\mathbb R^n$, $n\geq 3$:
{\em if the norm of the derivative $|D\varphi|$ belongs to Lorentz
space $L^{n-1,1}(D)$ and a mapping $\varphi: D\to D'$ has finite
distortion, then the inverse mapping belongs to Sobolev space
$W^1_{1,\loc}(D')$ and has finite distortion.} Recall that
$$
L^{n-1}(D)\subset L^{n-1,1}(D)\subset \bigcap\limits_{p>n-1}L^p(D).
$$

Note, that results about regularity of mappings inverse to Sobolev homeomorphisms follows from Theorem A. Indeed, substituting in the norm inequality for the inverse operator coordinate functions $x_j\in L^1_{p',\loc}(D)$ we see that $\varphi^{-1}$ belongs to $L^1_{1,\loc}(D')$.

The suggested method of investigation is based on a relation between Sobolev mappings,  composition operators of spaces of Lipschitz functions and a change of variable formula for weakly differentiable mappings.

\bigskip

\centerline{\bf 1.~Composition operators in Sobolev spaces}

\bigskip

A locally integrable function $f:D\to \mathbb R$ is {\it absolutely continuous on a
straight line $l$} having non-empty intersection with $D$ if it
is absolutely continuous on an arbitrary segment of this
 line which is contained in $D$. A function $f:D\to \mathbb
R$ belongs to the class $\ACL(D)$ ({\it absolutely continuous on
almost all straight lines}) if it is absolutely continuous on
almost all straight lines parallel to any coordinate axis.

Note that $f$ belongs to Sobolev space $L^1_{1}(D)$ if
and only if $f$ is locally integrable and it can be changed by a standard procedure on a set of measure zero (changed to its Lebesgue values at any point where the Lebesgue values exist) so , that a modified
function belongs to
$\ACL(D)$, and its partial derivatives $\frac{\partial f}{\partial
x_i}(x)$, $i=1,\ldots,n$, exist almost everywhere and are
integrable in $D$. From this point we will use such modified functions only.
Note that first weak derivatives of
the function $f$ coincide almost everywhere with the usual
partial derivatives (see, e.g., [14] ).

A mapping $\varphi:D\to \mathbb R^n$ belongs to the class $\ACL(D)$, if
its coordinate functions $\varphi_j$ belong to $\ACL(D)$, $j=1,\dots,n$.

We will use the notion of approximate differentiability.
Let $A$ be a subset of $\mathbb R^n$. Density of set $A$ at a point $x\in \mathbb R^n$ is the limit
$$
\lim\limits_{r\to 0}\frac{|B(x,r)\cap A|}{|B(x,r)|}.
$$
Here by symbol $|A|$ we denote Lebesgue measure of the set $A$.

A linear mapping $L: \mathbb R^n\to \mathbb R^n$ is called an approximate differential of a mapping $\varphi : D\to \mathbb R^n$
at point $a\in D$, if for every $\varepsilon>0$ the density of the set
$$
A_{\varepsilon}=\{x\in D : |\varphi(x)-\varphi(a)-L(x-a)|<\varepsilon |x-a| \}
$$
at point $a$ is equal to one.

A point $y\in \mathbb R^n$ is called an approximate limit of a mapping $\varphi : D\to \mathbb R^n$ at a point $x$, if the density of the set $D\setminus\varphi^{-1}(W)$ at this point is equal to zero for every neighborhood $W$ of the point $y$. 

For a mapping $\varphi : D\to\mathbb R^n$ we define approximate partial derivatives
$$
\ap\frac{\partial \varphi_i}{\partial x_j}(x)=\ap\lim\limits_{t \to 0}\frac{\varphi_i(x+te_j)-\varphi_i(x)}{t},\quad i,j=1,...,n.
$$

Approximate differentiable mappings are closely connected with Lipschitz mappings.
Recall, that a mapping $\varphi : D\to \mathbb R^n$ is a Lipschitz mapping, if there exists a constant $K<+\infty$ such that
$$
|\varphi(x)-\varphi(y)|\leq K |x-y|
$$
for every points $x,y \in D$.

The value 
$$
\|\varphi\mid\Lip (D)\| =\sup\limits_{x,y\in D} \frac{|\varphi(x)-\varphi(y)|}{|x-y|}
$$
we call the norm of $\varphi$ in the space $\Lip (D)$.

The next assertion describes this connection between approximate differentiable mappings and Lipschitz mappings in details [15].

\vskip 0.3cm
{\bf Theorem~1.}
Let $\varphi : D\to\mathbb R^n$ be a measurable mapping. Then the following assertions are equivalent:

\noindent
1) The mapping $\varphi : D\to\mathbb R^n$ is approximate differentiable almost everywhere in $D$.

\noindent
2) The mapping $\varphi : D\to\mathbb R^n$ has  approximate partial derivatives $\ap\frac{\partial \varphi_i}{\partial x_j}$, $i,j=1,...,n$ almost everywhere in $D$.

\noindent
3) There exists a collection of closed sets $\{A_k\}_{k=1}^{\infty}$, $A_k\subset A_{k+1}\subset D$, such that a restriction $\varphi \vert_{A_k}$ is a Lipschitz mapping on the set $A_k$ and 
$$
\biggl|D\setminus\sum\limits_{k=1}^{\infty}A_k\biggr|=0.
$$

\vskip 0.3cm

If a mapping $\varphi:D\to D^{\prime}$ has approximate partial
derivatives $\ap\frac{\partial \varphi_i}{\partial x_j}$ almost everywhere
in $D$, $i,j=1,\dots,n$, then  the formal Jacobi matrix
$D\varphi(x)=(\ap\frac{\partial \varphi_i}{\partial x_j}(x))$, $i,j=1,\dots,n$,
and its Jacobian determinant $J(x,\varphi)=\det D\varphi(x)$ are well defined at
almost all points of $D$. The norm $|D\varphi(x)|$ of the matrix
$D\varphi(x)$ is the norm of the linear operator determined by the
matrix in Euclidean space $\mathbb R^n$.

In the theory of mappings with bounded mean distortion additive set functions play a significant role.
Let us recall that a nonnegative mapping $\Phi$ defined on open subsets
of $D$ is called a {\it finitely
quasiadditive} set function [16] if

1) for any point $x\in D$, there exists $\delta$,
$0<\delta<\dist(x,\partial D)$, such that $0\leq
\Phi(B(x,\delta))<\infty$ (here and in what follows
$B(x,\delta)=\{y\in\mathbb R^n: |y-x|<\delta\}$);

2) for any  finite collection $U_i\subset U\subset D$,
$i=1,\dots,k$ of mutually disjoint open sets the following
inequality $\sum\limits_{i=1}^k \Phi(U_i)\leq \Phi(U)$ takes
place.

Obviously, the last inequality can be extended to a countable collection of mutually
disjoint open sets from $D$, so a finitely quasiadditive set
function is also {\it countable quasiadditive.}

If instead of the second condition we suppose that for any finite
collection $U_i\subset D$, $i=1,\dots,k$ of mutually disjoint
open subsets of $D$ the equality
$$
\sum\limits_{i=1}^k \Phi(U_i)= \Phi(U)
$$
takes place, then such set function is said to be {\it finitely
additive}. If the last equality can be extended to a
countable collection of mutually disjoint open subsets of $D$,
then such set function is said to be {\it countable additive.}

A nonnegative mapping $\Phi$ defined on open subsets of $D$ is called a {\it monotone} set function [16] if
$\Phi(U_1)\leq\Phi(U_2)$ under the condition, that $U_1\subset
U_2\subset D$ are open sets.

Note, that a monotone (countable) additive set function is the (countable) quasiadditive set function.

Let us reformulate an auxiliary result from [16] in a convenient for this study way.

\vskip 0.3cm

{\bf Proposition~1.} Let a monotone finitely additive set
function $\Phi$ be defined on open subsets of the domain
$D\subset\mathbb R^n$. Then for almost all points $x\in D$ the
volume derivative
$$
\Phi'(x)=\lim\limits_{\delta\to 0, B_{\delta}\ni x}
\frac{\Phi(B_{\delta})}{|B_{\delta}|}
$$
is finite and for any open set $U\subset D$, the inequality
$$
\int\limits_{U}\Phi'(x)~dx\leq \Phi(U)
$$
is valid.
\vskip 0.3cm

A nonnegative finite valued set function $\Phi$ defined on a collection of
measurable subsets of an open set $D$ is
said to be {\it absolutely continuous} if for every number
$\varepsilon>0$ can be found a number $\delta>0$ such that
$\Phi(A)<\varepsilon$ for any measurable sets $A\subset D$ from
the domain of definition of $\Phi$, which satisfies the condition
$|A|<\delta$.

Let $E$ be a measurable subset of $\mathbb R^n$, $n\geq 2$. Define Lebesgue space $L_p(E)$, $1\leq p\leq\infty$, as a Banach space of locally summable  functions $f:E\to \mathbb R$ equipped with the following norm:
$$
\|f\mid L_p(E)\|=
\biggr(\int\limits_E|f|^p(x)\,dx\biggr)^{1/p},\,\,\,1\leq p<\infty,
$$
and
$$
\|f\mid L_{\infty}(E)\|=
\ess\sup\limits_{x\in E}|f(x)|,\quad p=\infty.
$$
A function $f$ belongs to the space $L_{p,\loc}(E)$, $1\leq p\leq\infty$, if $f\in L_{p}(F)$ for every compact set $F\subset E$.

For an open subset $D\subset\mathbb R^n$ define the seminormed Sobolev space $L^1_p(D)$, $1\leq p\leq\infty$,
as a  space  of locally summable, weakly differentiable functions $f:D\to\mathbb R$ equipped with the following seminorm:
$$
\|f\mid L^1_{p}(D)\|=\| \nabla f\mid L_p(D)\|, \,\,\,\,1\leq p\leq\infty.  
$$
Here $\nabla f$ is the weak gradient of the function $f$, i.~e. $ \nabla f = (\frac{\partial f}{\partial x_1},...,\frac{\partial f}{\partial x_n})$,
$$
\int\limits_D f \frac{\partial \eta}{\partial x_i}~dx=-\int\limits_D \frac{\partial f}{\partial x_i} \eta~dx, \quad \forall \eta\in C_0^{\infty}(D),\quad i=1,...,n.
$$ 
As usual $C_0^{\infty}(D)$ is the space of infinitely smooth functions with a compact support.

Note, that smooth functions are dense in $L^1_p(D)$, $1\leq p<\infty$ (see, for example [14],\, [17]). If $p=\infty$ we can assert only that for arbitrary function $f\in L^1_p(D)$ there exists a sequence of smooth functions $\{f_k\}$ converges locally uniformly to $f$ and 
$\|f_k\mid L^1_{\infty}(D)\|\to \|f\mid L^1_{\infty}(D)\|$ (see [17]). 

The Sobolev space $W^1_p(D)$, $1\leq p\leq\infty$, is a Banach space of locally summable, weakly differentiable functions $f:D\to\mathbb R$,  equipped with the following norm:
$$
\|f\mid W^1_p(D)\|= \|f\mid L_p(D)\|+\|f\mid L^1_p(D)\|.
$$

A function $f$ belongs to the space $L^1_{p,\loc}(D)$ ($W^1_{p,\loc}(D)$), $1\leq p\leq\infty$, if $f\in L^1_{p}(K)$ ($f\in W^1_{p}(K)$) for every compact subset $K\subset D$.
The Sobolev space $\overset\circ{L}^1_p(D)$ is the
closure of the space $C^{\infty}_0(D)$ in ${L}^1_p(D)$.

A mapping $\varphi:D\to D^{\prime}$ belongs to Lebesgue class $L_p(E)$ if
its coordinate functions $\varphi_j$, $j=1,\dots,n$ belong to $L_p(E)$.
A mapping $\varphi:D\to D^{\prime}$ belongs to Sobolev class $W^1_p(D)$ ($L^1_p(D)$) if
its coordinate functions $\varphi_j$, $j=1,\dots,n$, belong to $W^1_p(D)$ ($L^1_p(D)$).

We say that a mapping $\varphi : D\to D'$ generates a bounded composition operator 
$$
\varphi^{\ast} : L^1_p(D')\to L^1_q(D), \,\,\,1\leq q\leq p\leq\infty,
$$
if for every function $f\in L^1_p(D')$ the composition $f\circ\varphi \in L^1_q(D)$ and the inequality
$$
\|\varphi^{\ast}f \mid L^1_q(D)\|\leq K \|f\mid L^1_p(D')\|
$$
holds.

{\bf Theorem~2.} A homeomorphism $\varphi : D\to D^{\prime}$
between two domains $D,D^{\prime} \subset \mathbb R^n $ generates a bounded composition operator
$$
\varphi^{\ast}: L^1_{\infty}(D')\to L^1_q(D),\,\, 1< q<+\infty,
$$
if and only if $\varphi$ belongs to the Sobolev space $L^1_q(D)$.
\vskip 0.3cm

{\sc Proof.} 
{\it Necessity.}
Substituting in the inequality
$$
\|\varphi^{\ast}f\mid L^1_q(D)\|\leq K \|f\mid L^1_{\infty}(D')\|
$$
the test functions $f_j(y)=y_j\in L^1_{\infty}(D')$, $j=1,...,n$ we see that $\varphi$ belongs to $L^1_q(D)$.

{\it Sufficiency.}
Let a function $f\in L^1_{\infty}(D')\cap C^{\infty}(D')$. Then 
\begin{multline}
\|\varphi^{\ast}f\mid L^1_q(D)\|=\biggl(\int\limits_D|\nabla (f\circ\varphi)|^q~dx\biggr)^{\frac{1}{q}}
\leq\biggl(\int\limits_D |D\varphi|^q|\nabla f|^q(\varphi(x))~dx\biggr)^{\frac{1}{q}}\\
\leq\biggl(\int\limits_D |D\varphi)|^q~dx\biggr)^{\frac{1}{q}}\|f\mid L^1_{\infty}(D')\|=
\|\varphi\mid L^1_q(D)\|\cdot \|f\mid L^1_{\infty}(D')\|.
\nonumber
\end{multline}
For arbitrary function $f\in L^1_{\infty}(D')$ consider a sequence of smooth functions $f_k\in L^1_{\infty}(D')$ such that
$$
\lim\limits_{k\to\infty}\|f_k\mid L^1_{\infty}(D')\|=\|f\mid L^1_{\infty}(D')\|
$$
and $f_k$ converges locally uniformly to $f$ in $D'$. Then, the sequence $\varphi^{\ast}f_k$ converges locally uniformly to $\varphi^{\ast}f$ in $D$ and
is a bounded sequence in $L^1_q(D)$. Since the space $L^1_q(D)$, $1<q<\infty$, is a reflexive space there exists a subsequence $f_{k_l}\in L^1_q(D)$ which weakly converges to $f\in L^1_q(D)$ and
$$
\|\varphi^{\ast}f\mid L^1_q(D)\|\leq \liminf\limits_{l\to\infty}\|\varphi^{\ast}f_{k_l}\mid L^1_q(D)\|.
$$
So, passing to limit when $l$ tends to $+\infty$ in the inequality
$$
\|\varphi^{\ast}f_{k_l}\mid L^1_q(D)\|\leq K \|f_{k_l}\mid L^1_{\infty}(D')\|
$$
we obtain
$$
\|\varphi^{\ast}f\mid L^1_q(D)\|\leq K \|f\mid L^1_{\infty}(D')\|.
$$
\vskip 0.3cm

The next theorem gives a "localization" property of the composition operator on spaces of functions with compact support and/or its closure in $L^1_{\infty}$.

\vskip 0.3cm

{\bf Theorem~3.} Let a homeomorphism $\varphi : D\to D^{\prime}$
between two domains $D,D^{\prime} \subset \mathbb R^n $ generates a bounded composition operator
$$
\varphi^{\ast}: L^1_{\infty}(D')\to L^1_q(D),\,\, 1\leq q<+\infty.
$$
Then there exists a bounded
monotone countable additive function $\Phi(A')$ defined on open bounded
subsets of $D'$ such that for every function $f\in
\overset{\circ}{L}^1_{\infty}(A')$ the inequality
$$
\int\limits_{\varphi^{-1}(A)}|\nabla(f\circ \varphi)|^q~dx\leq
\Phi(A'){\esssup\limits_{y\in A'}|\nabla f|^q(y)}
$$
holds.
\vskip 0.3cm

{\sc Proof.} Let us define $\Phi(A')$ by the following way [2,\,6]
$$
\Phi(A^{\prime})=\sup\limits_{f\in
\overset{\circ}{L}_{\infty}^{1}(A^{\prime})} \Biggl(
\frac{\bigl\|\varphi^{\ast} f\mid {L}_{q}^{1}(D)\bigr\|}
{\bigl\|f\mid
\overset{\circ}{L}_{\infty}^{1}(A^{\prime})\bigr\|}
\Biggr)^{q},
$$

Let $A_1^{\prime}\subset A_2^{\prime}$ be bounded open
subsets of $D^{\prime}$. Extending
functions of space
$\overset{\circ}{L}_{\infty}^{1}(A_1^{\prime})$ by zero onto the set
$A_2^{\prime}$, we obtain an inclusion 
$\overset{\circ}{L}_{\infty}^{1}(A_1^{\prime})\subset
\overset{\circ}{L}_{\infty}^{1}(A_2^{\prime})$. Obviously 
$$
\|f\mid \overset{\circ}{L}_{\infty}^{1}(A_1^{\prime})\|=\|f\mid \overset{\circ}{L}_{\infty}^{1}(A_2^{\prime})\|
$$
for every $f\in \overset{\circ}{L}_{\infty}^{1}(A_1^{\prime})$. By  the following
inequality
\begin{multline}
\Phi(A_1^{\prime})=\sup\limits_{f\in
\overset{\circ}{L}_{\infty}^{1}(A_1^{\prime})} \Biggl(
\frac{\bigl\|\varphi^{\ast} f\mid {L}_{q}^{1}(D)\bigr\|}
{\bigl\|f\mid
\overset{\circ}{L}_{\infty}^{1}(A_1^{\prime})\bigr\|}
\Biggr)^{q}
=
\sup\limits_{f\in
\overset{\circ}{L}_{\infty}^{1}(A_1^{\prime})} \Biggl(
\frac{\bigl\|\varphi^{\ast} f\mid {L}_{q}^{1}(D)\bigr\|}
{\bigl\|f\mid
\overset{\circ}{L}_{\infty}^{1}(A_2^{\prime})\bigr\|}
\Biggr)^{q}\\
\leq\sup\limits_{f\in
\overset{\circ}{L}_{\infty}^{1}(A_2^{\prime})} \Biggl(
\frac{\bigl\|\varphi^{\ast} f\mid {L}_{q}^{1}(D)\bigr\|}
{\bigl\|f\mid
\overset{\circ}{L}_{\infty}^{1}(A_2^{\prime})\bigr\|}
\Biggr)^{q}=\Phi(A_2^{\prime}).
\nonumber
\end{multline}
 the set function $\Phi$ is monotone.

Let $A_{i}^{\prime}$, $i\in \mathbb N$, be open disjoint subsets at the domain
$D^{\prime}$,
$A^{\prime}_0=\bigcup\limits_{i=1}^{\infty}A_i^{\prime}$.
Choose arbitrary functions $f_i\in\overset{\circ}{L}_{\infty}^{1}(A_i^{\prime})$ with following
properties
$$
\bigl\|\varphi^{\ast} f_i\mid{L}_{q}^{1}
(D)\bigr\|\geq \bigl(\Phi(A_i^{\prime})\bigl(1-\frac{\varepsilon}{2^i}\bigr)\bigr)^{\frac{1}{q}}
\bigl\|f_i\mid\overset{\circ}{L}_{\infty}^{1}(A_i^{\prime})
\bigr\|
$$ 
and
$$
\bigl\|f_i\mid\overset{\circ}{L}_{\infty}^{1}(A_i^{\prime})
\bigr\|=1,
$$
while $i\in\mathbb N$. Here $\varepsilon\in(0,1)$ is a fixed number.
Letting $g_N=\sum\limits_{i=1}^{N}f_i$ we obtain
\begin{multline}
\bigl\|\varphi^{\ast} g_N\mid{L}_{q}^{1}
(D)\bigr\|\geq \biggl(\sum\limits_{i=1}^{N}
\left(\Phi(A_i^{\prime})\left(1-\frac{\varepsilon}{2^i}\right)\right)
\bigl\|f_i\mid\overset{\circ}{L}_{\infty}^{1}(A_i^{\prime})
\bigr\|^q\biggr)^{1/q}\\
= \biggl(\sum\limits_{i=1}^{N}\Phi(A_i^{\prime})
\left(1-\frac{\varepsilon}{2^i}\right)\biggr)^{\frac{1}{q}}
\biggl\|g_N\mid\overset{\circ}{L}_{\infty}^{1}
\Bigl(\bigcup\limits_{i=1}^{N}A_i^{\prime}\Bigr)\biggr\|
\\
\geq \biggl(\sum\limits_{i=1}^{N}\Phi(A_i^{\prime})
-\varepsilon\Phi(A_0^{\prime}) \biggr)^{\frac{1}{q}}
\biggl\|g_N\mid\overset{\circ}{L}_{\infty}^{1}
\Bigl(\bigcup\limits_{i=1}^{N}A_i^{\prime}\Bigr)\biggr\|
\nonumber
\end{multline}
since sets, on which the gradients $\nabla\varphi^{\ast} f_i$
do not vanish, are disjoint. From the last inequality  follows that
$$
\Phi(A_0^{\prime})^{\frac{1}{q}}\geq\sup\frac
{\bigl\|\varphi^{\ast} g_N\mid{L}_{q}^{1} (D)\bigr\|}
{\biggl\|g_N\mid\overset{\circ}{L}_{\infty}^{1}
\Bigl(\bigcup\limits_{i=1}^{N}A_i^{\prime}\Bigr)\biggr\|}\geq
\biggl(\sum\limits_{i=1}^{N}\Phi(A_i^{\prime})-\varepsilon\Phi(A_0^{\prime})
\biggr)^{\frac{1}{q}}.
$$
Here the upper bound is taken over all above-mentioned
functions $$
g_N\in\overset{\circ}{L}_{\infty}^{1}
\Bigl(\bigcup\limits_{i=1}^{N}A_i^{\prime}\Bigr).
$$
Since both
$N$ and $\varepsilon$ are arbitrary, we have finally
$$
\sum\limits_{i=1}^{\infty}\Phi(A^{\prime}_i) \leq
\Phi\Bigl(\bigcup\limits_{i=1}^{\infty}A^{\prime}_i\Bigr).
$$

The validity of the inverse inequality can be proved in a straightforward manner.
Indeed, choose functions $f_i\in\overset{\circ}{L}_{\infty}^{1}(A_i^{\prime})$ such that $\bigl\|f_i\mid\overset{\circ}{L}_{\infty}^{1}(A_i^{\prime})
\bigr\|=1$.

Letting $g=\sum\limits_{i=1}^{\infty}f_i$ we obtain
$$
\bigl\|\varphi^{\ast} g\mid{L}_{q}^{1}(D)\bigr\|\leq \biggl(\sum\limits_{i=1}^{\infty}
\Phi(A_i^{\prime})
\bigl\|f_i\mid\overset{\circ}{L}_{\infty}^{1}(A_i^{\prime})
\bigr\|^q\biggr)^{1/q}
= \biggl(\sum\limits_{i=1}^{\infty}\Phi(A_i^{\prime})
\biggr)^{\frac{1}{q}}\biggl\|g_N\mid\overset{\circ}{L}_{\infty}^{1}
\Bigl(\bigcup\limits_{i=1}^{\infty}A_i^{\prime}\Bigr)\biggr\|,
$$
since sets, on which the gradients $\nabla\varphi^{\ast} f_i$
do not vanish, are disjoint. From this inequality follows that
$$
\Phi\biggl(\bigcup\limits_{i=1}^{\infty}A_i^{\prime}\biggr)^{\frac{1}{q}}\leq\sup\frac
{\bigl\|\varphi^{\ast} g\mid{L}_{q}^{1} (D)\bigr\|}
{\biggl\|g\mid\overset{\circ}{L}_{\infty}^{1}
\Bigl(\bigcup\limits_{i=1}^{\infty}A_i^{\prime}\Bigr)\biggr\|}\leq
\biggl(\sum\limits_{i=1}^{\infty}\Phi(A_i^{\prime})
\biggr)^{\frac{1}{q}},
$$
where the upper bound is taken over all functions $g\in\overset{\circ}{L}_{\infty}^{1}
\Bigl(\bigcup\limits_{i=1}^{\infty}A_i^{\prime}\Bigr)$.

By the definition of the set function $\Phi$ we have
$$
\|\varphi^{\ast} f\mid L^1_q(D)\|^p\leq \Phi(A') \|f\mid\overset{\circ}{L}_{\infty}^{1}(A')\|^q
$$
Since the support of the function $f\circ\varphi$ is contained in the set $\varphi^{-1}(A')$ we have
$$
\int\limits_{\varphi^{-1}(A)}|\nabla(f\circ \varphi)|^q~dx\leq
\Phi(A'){\esssup\limits_{y\in A'}|\nabla f|^q(y)}.
$$

Theorem proved.

\vskip0.3cm

We recall some basic facts about $p$-capacity. Let $G\subset\mathbb
R^n$ be an open set and $E\subset G$ be a compact set. For $1\leq
p\leq\infty$ the $p$-capacity of the ring $(E,G)$ is defined as
$$
\cp_p(E,G)=\inf\bigl\{\int_G |\nabla u|^p : u\in
{L}^1_p(G)\cap C^{\infty}_0(G),\, u\geq 1\,\, \text{on}
\,\,E\bigr\}.
$$
Functions $u\in {L}^1_p(G)\cap C^{\infty}_0(G),\, u\geq 1\,\, \text{on}
\,\,E$, are called admissible functions for ring $(E,G)$.

We need the following estimate of the $p$-capacity [18].

\vskip 0.3cm
{\bf Lemma~1.} Let $E$ be a connected closed subset of
an open bounded set $G\subset\mathbb R^n, n\geq 2$, and $n-1<p<\infty$. Then
$$
\cp_p^{n-1}(E,G)\geq c\frac{(\diam E)^p}{|G|^{p-n+1}},
$$
where a constant $c$ depends on $n$ and $p$ only.

\vskip 0.3cm

For readers convenience we will prove this fact.
 
{\sc Proof.} Let $d$ be diameter of set $E$. Without loss of generality we can suggest, that $d=\dist(0, a)$ for some point $a=(0,..,0, a_n)$. For arbitrary number $t$, $0<t<d$, denote by $P_t$ the hyperplane $x_n=t$.

In the subspace $x_n=0$ we consider the unit $(n-2)$-dimensional sphere $S^{n-2}$ with the center at the origin
and fix an arbitrary point $z\in E\cap P_t$. For every point $y\in S^{n-2}$ denote by $R(y)$ the supremum of numbers $r_0$ such that $z+ry\in G$ while $0\leq r\leq r_0$.
Then for every admissible function $f\in C^{\infty}_0(G)$ the following inequality
$$
1= f(z)-f(z+R(y)y)\leq  \int\limits_0^{R(y)}|\nabla f(z+ry)|~dr=\int\limits_0^{R(y)}(|\nabla f(z+ry)|r^{\frac{n-2}{p}})r^{-\frac{n-2}{p}}~dr
$$
holds.
Applying H\"older inequality to the right side of the last inequality, we have
$$
1\leq \biggl(\frac{p-1}{p-n+1}\biggr)^{p-1}\bigl(R(y)\bigr)^{p-n+1}\int\limits_0^{R(y)}|\nabla f(z+ry)|^p r^{n-2}~dr.
$$

Multiplying both sides of this inequality on 
$((p-1)/(p-n+1))^{1-p}\cdot (R(y))^{n-p-1}$ 
and integrating by $y\in S^{n-2}$, we obtain
\begin{multline}
\biggl(\frac{p-1}{p-n+1}\biggr)^{p-1}\int\limits_{S^{n-2}} \bigl(R(y)\bigr)^{p-n+1}~dy\\
\leq\int\limits_{S^{n-2}}~dy
\int\limits_0^{R(y)}|\nabla f(z+ry)|^p r^{n-2}~dr\leq\int\limits_{P_t}|\nabla f|^p~dz.
\nonumber
\end{multline}

For the lower estimate of the left integral we use again H\"older inequality. Denote by $\omega_{n-2}$ the $n-2$-dimensional area of sphere $S^{n-2}$. By simple calculations we get
\begin{multline}
\omega_{n-2}^p=\biggl(\int\limits_{S^{n-2}}~dy\biggr)^p\leq \biggl(\int\limits_{S^{n-2}}\bigl(R(y)\bigr)^{n-p-1}~dy\biggr)^{n-1}
\biggl(\int\limits_{S^{n-2}}\bigl(R(y)\bigr)^{n-1}~dy\biggr)^{p+1-n}\\
\leq((n-1)m_{n-1}(G\cap P_t))^{p-n+1}\biggl(\int\limits_{S^{n-2}}\bigl(R(y)\bigr)^{n-p-1}~dy\biggr)^{n-1}.
\nonumber
\end{multline}
Here $m_{n-1}(A)$ is $(n-1)$-Lebesgue measure of the set $A$.

Denote by $u(t)=m_{n-1}(G\cap P_t)$. Using the last estimate we obtain
$$
\int\limits_{P_t} |\nabla f|^p~dz\geq \biggr( \frac{p-1}{p-n+1}\biggl)^{1-p}(n-1)^{\frac{n-p-1}{n-1}}\omega_{n-2}^{\frac{p}{n-1}}\bigl(u(t)\bigr)^{\frac{n-p-1}{n-1}}.
$$
After integrating by $t\in (0,d)$ we have
$$
\int\limits_{G} |\nabla f|^p~dx\geq \biggr( \frac{p-1}{p-n+1}\biggl)^{1-p}(n-1)^{\frac{n-p-1}{n-1}}\omega_{n-2}^{\frac{p}{n-1}}\int\limits_0^d\bigl(u(t)\bigr)^{\frac{n-p-1}{n-1}}~dt.
$$
By H\"older inequality
\begin{multline}
d^p=\biggl(\int\limits_0^d~dt\biggr)^p\leq\biggl(\int\limits_0^du(t)~dt\biggr)^{p-n+1} 
\biggl(\int\limits_0^d\bigl(u(t)\bigr)^{\frac{n-p-1}{n-1}}~dt\biggr)^{n-1} \\
\leq |G|^{p-n+1}\biggl(\int\limits_0^d\bigl(u(t)\bigr)^{\frac{n-p-1}{n-1}}~dt\biggr)^{n-1}.
\nonumber
\end{multline}

Therefore
$$
\int\limits_{G}|\nabla f|^p~dx\geq \biggr( \frac{p-1}{p-n+1}\biggl)^{1-p}(n-1)^{\frac{n-p-1}{n-1}}\omega_{n-2}^{\frac{p}{n-1}}\biggl(\frac{d^p}{|G|^{p-n+1}}\biggr)^{\frac{1}{n-1}}.
$$
Since $f$ is an arbitrary admissible function the required inequality is proved.

\vskip 0.3cm

Let us define a class $BVL$ of mappings with finite variation.
A mapping $\varphi : D\to\mathbb R^n$ belongs to the class $\BVL(D)$
(i.e., has {\it finite variation on almost all straight lines})
if it has finite variation on almost all straight lines $l$ parallel
to any coordinate axis: for any finite number of points $t_1,...,t_k$ that belongs to such straight line $l$
$$
\sum\limits_{i=0}^{k-1}|\varphi(t_{i+1})-\varphi(t_i)|<+\infty.
$$

For a mapping $\varphi$ with finite variation on almost all straight lines, the partial
derivatives $\partial \varphi_i/\partial x_j$, $i,j=1,\dots,n$,
exists almost everywhere in $D$.

\vskip 0.3cm
{\bf Theorem~4. [11]} Let a homeomorphism $\varphi : D\to D^{\prime}$
generates a bounded composition operator
$$
\varphi^{\ast}: L^1_{\infty}(D')\to L^1_{q}(D),\,\,\,q>n-1.
$$
Then the inverse
homeomorphism $\varphi^{-1}:D'\to D$ belongs to the class $\BVL(D')$.
\vskip 0.3cm

For readers convenience we reproduce here a slightly modified proof of this fact.

{\sc Proof.} Take an arbitrary $n$-dimensional open
parallelepiped $P$ such that $\overline{P}\subset D'$ and its
edges are parallel to coordinate axis. Let us show that
$\varphi^{-1}$ has finite variation on almost all intersection of $P$
and straight lines parallel to $x_n$-axis.

Let $P_0$ be the projection of $P$ on the subspace $x_n=0$, and
let $I$ be the projection of $P$ on the coordinate axis $x_n$.
Then $P=P_0\times I$. The monotone countable-additive function
$\Phi$ determines a monotone countable additive function of open
sets $A\subset P_0$ by the rule $\Phi(A,P_0)=\Phi(A\times I)$. For
almost all points $z\in P_0$, the quantity
$$
\overline{\Phi'}(z,P_0) =\overline{\lim_{r\to
0}}\biggl[\frac{\Phi(B^{n-1}(z,r),P_0)} {r^{n-1}}\biggr]
$$
is finite [19] (here $B^{n-1}(z,r)$ is the $(n-1)$-dimensional ball of
radius $r>0$ centered at the point $z$).

The $n$-dimensional Lebesgue measure $\Psi(U)=|\varphi^{-1}(U)|$, where
$U$ is an open sen in $D'$, is a monotone countable additive
function and, therefore, also determines a monotone countable
additive function $\Psi(A,P_0) = \Psi(A\times I)$ defined on open
sets $A\subset P_0$. Hence  $\overline{\Psi'}(z,P_0)$ is finite for
almost all points $z\in P_0$.

Choose an arbitrary point $z\in P_0$ where 
 $\overline{\Phi'}(z,P_0)<+\infty$ and $\overline{\Psi'}(z,P_0)<+\infty$. On the section $I_z=\{z\}\times I$ of the
parallelepiped $P$, take arbitrary mutually disjoint closed
intervals $\Delta_1,...,\Delta_k$ with lengths $b_1,...,b_k$
respectively. Let $R_i$ denote the open set of points for which
distances from $\Delta_i$ smaller than a given $r>0$:
$$
R_i=\{x\in G : \dist(x,\Delta_i)<r\}.
$$
Consider
the ring $(\Delta_i, R_i)$. Let $r>0$ be selected so that
$r<cb_i$ for $i=1,\dots,k$, where $c$ is a sufficiently small constant.
Then the function $u_i(x)=\dist(x,\Delta_i)/r$ is an admissible
for ring $(\Delta_i, R_i)$. 

By Theorem~3 we have the estimate
$$
\|\varphi^{\ast} u_i \mid L^1_q(D)\|^q\leq \Phi(A') \|u_i \mid\overset{\circ}{L}_{\infty}^{1}(A')\|^q
$$
for every function $u_i$, $i=1,...,k$.

Hence, for every ring $(\Delta_i, R_i)$, $i=1,...,k$, the inequality
$$
\cp_q^{\frac{1}{q}}(\varphi^{-1}(\Delta_i), \varphi^{-1}(R_i))\leq \Phi(R_i)^{\frac{1}{q}}\cp_{\infty}(\Delta_i, R_i)
$$
holds. 

The function $u_i(x)=\dist(x,\Delta_i)/r$ is admissible
for ring $(\Delta_i, R_i)$ and we have the upper estimate
$$
\cp_{\infty}(\Delta_i, R_i)\leq |\nabla u_i|=\frac{1}{r}.
$$

Applying the
lower bound for the capacity of the ring (Lemma~1), we obtain
$$
\biggl(\frac{(\diam
\varphi^{-1}(\Delta_i))^{q/(n-1)}}{|\varphi^{-1}(R_i)|^{(q-n+1)/(n-1)}}\biggr)^{\frac{1}{q}}\leq
c_1 \Phi(R_i)^{\frac{1}{q}}\cdot\frac{1}{r}.
$$
This inequality gives
$$
\diam \varphi^{-1}(\Delta_i)\leq
c_2\biggl(\frac{|\varphi^{-1}(R_i)|}{r^{n-1}}\biggr)^{\frac{q-n+1}{q}}
\cdot \biggl(\frac{\Phi(R_i)}{r^{n-1}}\biggr)^{\frac{n-1}{q}}.
$$
Summing over $i=1,\dots,k$
we obtain
$$
\sum\limits_{i=1}^k\diam \varphi^{-1}(\Delta_i)\leq
c_2 \sum\limits_{i=1}^k \biggl(\frac{|\varphi^{-1}(R_i)|}{r^{n-1}}\biggr)^{\frac{q-n+1}{q}}
\cdot \biggl(\frac{\Phi(R_i)}{r^{n-1}}\biggr)^{\frac{n-1}{q}}.
$$

Hence
$$
\sum\limits_{i=1}^k\diam \varphi^{-1}(\Delta_i)\leq
c_2 \biggl(\sum\limits_{i=1}^k\frac{|\varphi^{-1}(R_i)|}{r^{n-1}}\biggr)^{\frac{q-n+1}{q}}
\cdot \biggl(\sum\limits_{i=1}^k\frac{\Phi(R_i)}{r^{n-1}}\biggr)^{\frac{n-1}{q}}.
$$

Using the Besicovitch type theorem [20] for the estimate of the
value of the function $\Phi$ in terms of the multiplicity of a
cover, we obtain
$$
\sum\limits_{i=1}^k\diam \varphi^{-1}(\Delta_i)\leq
c_3 \biggl(\frac{|\varphi^{-1}(\bigcup_{i-1}^k R_i)|}{r^{n-1}}\biggr)^{\frac{q-n+1}{q}}
\cdot \biggl(\frac{\Phi(\bigcup_{i-1}^k R_i)}{r^{n-1}}\biggr)^{\frac{n-1}{q}}.
$$
Hence
$$
\sum\limits_{i=1}^k \diam \varphi^{-1}(\Delta_i)\leq
c_3\biggl(\frac{|\varphi^{-1}(B^{n-1}(z,r),P_0)|}{r^{n-1}}\biggr)^{\frac{q-n+1}{q}}
\cdot
\biggl(\frac{\Phi(B^{n-1}(z,r),P_0)}{r^{n-1}}\biggr)^{\frac{n-1}{q}}.
$$
Because  $\overline{\Phi'}(z,P_0)<+\infty$ and $\overline{\Psi'}(z,P_0)<+\infty$ we obtain finally 

$$
\sum\limits_{i=1}^k \diam \varphi^{-1}(\Delta_i)<+\infty.
$$
Therefore $\varphi^{-1}\in \BVL(D')$.

Theorem proved.

\bigskip

\centerline{\bf 2.~Invertibility of composition operators}

\bigskip

Let us recall the change of variable formula for Lebesgue integral [21].
Let a mapping $\varphi : D\to \mathbb R^n$ be such that
there exists a collection of closed sets $\{A_k\}_1^{\infty}$, $A_k\subset A_{k+1}\subset D$ for which restrictions $\varphi \vert_{A_k}$ are Lipschitz mapping on sets $A_k$ and 
$$
\biggl|D\setminus\sum\limits_{k=1}^{\infty}A_k\biggr|=0.
$$
Then there exists a measurable set $S\subset D$, $|S|=0$ such that  the mapping $\varphi:D\setminus S \to \mathbb R^n$ has Luzin $N$-property and the change of variable formula
$$
\int\limits_E f\circ\varphi (x) |J(x,\varphi)|~dx=\int\limits_{\mathbb R^n\setminus \varphi(S)} f(y)N_f(E,y)~dy
$$ 
holds for every measurable set $E\subset D$ and every nonnegative Borel measurable function $f: \mathbb R^n\to\mathbb R$. Here 
$N_f(y,E)$ is the multiplicity function defined as the number of preimages of $y$ under $f$ in $E$.

If a mapping $\varphi$ possesses Luzin $N$-property (the image of a set of measure zero has measure zero), then $|\varphi (S)|=0$ and the second integral can be rewritten as the integral on $\mathbb R^n$.
Note, that if a homeomorphism $\varphi : D\to D'$ belongs to the Sobolev space $W^1_{n,\loc}(D)$ then $\varphi$ has Luzin $N$-property and the change of variable formula holds [22].

If a mapping $\varphi : D\to \mathbb R^n$ belongs to the Sobolev space $W^1_{1,\loc}(D)$ then by [21] there exists a collection of closed sets $\{A_k\}_1^{\infty}$, $A_k\subset A_{k+1}\subset D$ for which restrictions $f\vert_{A_k}$ are Lipschitz mapping on sets $A_k$ and 
$$
\biggl|D\setminus\sum\limits_{k=1}^{\infty}A_k\biggr|=0.
$$ 
Hence for such mappings the previous change of variable formula is correct.

Like in [23] (see also [13]) we define a measurable function

$$
\mu(y)=
\begin{cases}
\biggl(\frac{|\adj D\varphi|(x)}{|J(x,\varphi)|}\biggr)_{x=\varphi^{-1}(y)}
\quad\text{if}\quad x\in D\setminus S \quad\text{and}\quad J(x,\varphi)\ne 0,
\\
\,\,0\quad\quad\quad\quad\quad\quad \quad\quad \quad    \text{otherwise}.
\end{cases}
$$

Because the homeomorphism $\varphi$ has finite distortion the function $\mu(y)$ is well defined almost everywhere in $D'$.

The following lemma was proved (but does not formulated) in [13] under an additional assumption that $|D\varphi|$ belongs to the Lorentz space $L^{n-1,n}(D)$.

\vskip 0.3cm
{\bf Lemma~2.} Let a homeomorphism $\varphi:D\to D', \varphi(D)=D'$ belongs to the
Sobolev space $L^1_q(D)$ for some $q> n-1$. Then the function $\mu$ is locally
integrable in the domain $D'$.
\vskip 0.3cm

{\sc Proof.} 
Using the change of variable formula for Lebesgue integral [21] and Luzin $N$-property of $\varphi$  we have the following equality
$$
\int\limits_{D'}\mu(y)~dy=\int\limits_{D'\setminus \varphi(S)}\mu(y)~dy=\int\limits_{D\setminus S} | \mu(\varphi(x))|J(x,\varphi)|~dx= \int\limits_{D} |\adj D\varphi|(x)~dx.
$$
Applying H\"older inequality, we obtain that for every compact subset $F'\subset D'$
$$
\int\limits_{F'}\mu(y)~dy\leq \int\limits_{F} |\adj D\varphi|(x)~dx\leq C \int\limits_{F}
|D\varphi|^{n-1}(x)~dx,
$$
where $F'=\varphi(F)$.
Therefore,   $\mu$ belongs to $L_{1,\loc}(D')$, since $\varphi$ belongs to $L^1_q(D)$, $q>n-1$, and as consequence $\varphi\in L^1_{n-1,\loc}(D)$.
\vskip 0.3cm

{\bf Theorem~5.} Let a homeomorphism $\varphi : D\to D^{\prime}$, $\varphi(D)=D^{\prime}$, has finite distortion, Luzin $N$-property (the image of a set measure zero is a set measure zero) and  generates a bounded composition operator
$$
\varphi^{\ast}: L^1_{\infty}(D')\to L^1_q(D),\,\,\,q>n-1.
$$
Then the inverse
homeomorphism $\varphi^{-1}:D'\to D$  has integrable first weak derivatives and induces a bounded composition operator
$$
(\varphi^{-1})^{\ast} : L^1_{q'}(D)\to L^1_1(D'),\quad q'=q/(q-n+1).
$$
\vskip 0.3cm

{\sc Proof.} We prove that $\varphi^{-1}\in \ACL(D')$. Since absolute continuity is the local property, it is sufficient to prove that the mapping $\varphi^{-1}$ belongs to $\ACL$ on every compact subset of $D'$. Consider arbitrary cube $Q'\in D'$, $\overline{Q'}\in D$, with edges parallel to coordinate axes, and $Q=\varphi^{-1}(Q')$. For $i=1,\dots n$ we will use a notation: $Y_i=(x_1,...,x_{i-1},x_{i+1},...,x_n)$,
$$
F_i(x)=(\varphi_1(x),\dots,\varphi_{i-1}(x),\varphi_{i+1}(x),\dots,\varphi_n(x))
$$
and  $Q'_i$  is the intersection of the cube $Q'$ with a line
$Y_i=\const$. 

Using the change of variable formula and the Fubini theorem [24] we obtain the following estimate
$$
\int\limits_{F_i(Q)}H^{n-1}(dY_i)\int\limits_{Q_i'}\mu(y)~H^1(dy)=\int\limits_{Q'} \mu(y)~dy=\int\limits_Q |\adj D\varphi| (x)~dx<+\infty.
$$
Hence for almost all $Y_i\in F_i(Q)$
$$
\int\limits_{Q_i'}\mu(y)~H^1(dy)<+\infty.
$$

Let $\ap J\varphi(x)$ be an approximate Jacobian of
the trace of the mapping $\varphi$ on the set $\varphi^{-1}(Q'_i)$ [24].  Consider a point $x\in Q$ in which there exists a non-generated approximate differential $\ap Df(x)$ of the mapping $\varphi: D\to D'$. Let $L: \mathbb R^n\to \mathbb R^n$ be a linear mapping induced by this approximate differential $\ap Df(x)$. We denote by the symbol $P$ the image of the unit cube $Q_0$ under the linear mapping $L$ and by $P_i$ the intersection of $P$ with the image of the line $x_i=0$. Let $d_i$ be a length of $P_i$. Then 
$$
d_i\cdot|\adj DF_i|(x)=|Q_0|=|J(x,\varphi)|.
$$
So, since $d_i=\ap J\varphi(x)$ we obtain that for
almost all $x\in Q\setminus Z$, $Z=\{x\in D : J(x,\varphi)=0\}$, we have
$$
\ap J\varphi(x)=\frac{|J(x,\varphi)|}{|\adj DF_i|(x)}.
$$

So, we have for arbitrary
compact set $A'\subset Q'_i$, and for almost all $Y_i\subset F_i(Q)$,
the following inequality:
\begin{multline}
 H^1(\varphi^{-1}(A'))\leq \int\limits_{\varphi^{-1}(A')}\frac{|\adj
D\varphi|(x)}{|\adj
DF_i|(x)}~H^1(dx)\\
=\int\limits_{\varphi^{-1}(A')}\frac{|\adj
D\varphi|(x)}{|J(x,\varphi)|}\cdot\frac{|J(x,\varphi)|}{|\adj
DF_i|(x)}~H^1(dx)
=
\int\limits_{\varphi^{-1}(A')}\mu(\varphi(x))\ap
J\varphi(x)~H^1(dx).
\nonumber
\end{multline}

By using the change of variable formula for the Lebesgue integral [24, 25] we obtain
$$
H^1(f^{-1}(A'))\leq \int\limits_{A'}\mu(y)~H^1(dy)<+\infty.
$$
Therefore, the mapping $\varphi^{-1}$ is absolutely continuous on almost all
lines in $D'$and is a weakly differentiable mapping.

Since the homeomorphism $\varphi$ has Luzin $N$-property then preimage of a set positive measure is a set positive measure. Hence, the volume derivative of the inverse mapping
$$
J_{\varphi^{-1}}(y)=\lim\limits_{r\to 0}\frac{|\varphi^{-1}(B(y,r))|}{|B(y,r)|}>0
$$
almost everywhere in $D'$.
So $J(y,\varphi^{-1}) \neq 0$ for almost all points $y\in D$. Integrability of the $q'$-distortion follows from the inequality
$$
|D\varphi^{-1}|(y)\leq |D\varphi(x)|^{n-1}\big/ |J(x,\varphi)|
$$
which holds for almost all points $y=\varphi(x)\in D'$.

Indeed, with the help of the change of variable formula, we have
\begin{multline}
\int\limits_{D'}\biggl(\frac{|D\varphi^{-1}(y)|^{q'}}{|J(y,\varphi^{-1})|}\biggr)^{\frac{1}{q'-1}}~dy=
\int\limits_{D'}\biggl(\frac{|D\varphi^{-1}(y)|}{|J(y,\varphi^{-1})|}\biggr)^{\frac{q'}{q'-1}}|J(y,\varphi^{-1})|~dy\\
\leq \int\limits_{D}\biggl(\frac{|D\varphi^{-1}(\varphi(x))|}{|J(\varphi(x),\varphi^{-1})|}\biggr)^{\frac{q'}{q'-1}}~dx\leq \int\limits_{D}|D\varphi(x)|^q~dx<+\infty,
\nonumber
\end{multline}
since by Theorem~2 $\varphi$ belongs to $L^1_q(D)$.

The boundedness of the composition operator follows from integrability of the $p'$-distortion [2].
The theorem proved.

\vskip 0.5cm
We are pleasure to thank Professor Jan Maly for helpful discussions.

\vskip 0.3cm

\centerline{REFERENCES}

\begin{enumerate}

\item
Vodop'yanov S.~K. and Gol'dshtein V.~M. {\it Structure isomorphisms of spaces $W^1_n$ and quasiconformal mappings.}//
Siberian Math. J. -- 1975. -- V.~16. -- P.~224--246.

\item
Ukhlov A.~D. {\it Mappings that generate embeddings of Sobolev spaces.}//
Siberian Math. J. -- 1993. -- V.~34. -- N.~1. -- P.~165--171.

\item
Gol'dshtein V., Gurov L. {\it Applications of change of variable operators for exact embedding theorems.}//
Integral equations operator theory -- 1994. -- V.~19. -- N.~1. -- P.~1--24.

\item
Ball J.~M. {\it Convexity conditions and existence theorems in nonlinear elasticity.}//
Arch. Rat. Mech. Anal. -- 1976. -- V.~63. -- P.~337--403.

\item
Ball J.~M. {\it Global invertability of Sobolev functions and the interpretation of matter.}//
Proc. Roy. Soc. Edinburgh -- 1981. -- V.~88A. --  P.~315--328.

\item
Vodop'yanov S.~K., Ukhlov A.~D. {\it Sobolev spaces and $(p,q)$-quasiconformal mappings of Carnot groups.}//
Siberian Math. J. -- 1998. -- V.~39. -- N.~4. -- P.~776--795.

\item
Vodop'yanov S.~K., Ukhlov A.~D. {\it Mappings with bounded $(p,q)$-distortion on Carnot groups.}//
Bull. Sci. Math. (to appear)

\item
Gol'dshtein V., Ramm A.~G. {\it Compactness of the embedding operators for rough domains.}//
Math. Inequalities and Applications -- 2001. -- V.~4. -- N.~1. -- P.~127--141.

\item
Gol'dshtein V., Ukhlov A. {\it Weighted Sobolev spaces and embedding theorems.}//
Transactions of Amer. Math. Soc. (to appear)

\item
Muller S., Tang Q., Yan B.~S. {\it On a new class of elastic deformations not allowing for cavitation.}//
Ann. Inst. H. Poincare. Anal. non. lineaire -- 1994. -- V.~11. -- N.~2. -- P.~217--243.

\item
Ukhlov~A. {\it Differential and geometrical properties of Sobolev mappings.}//
Mathematical Notes -- 2004. -- V.~75. -- N.~2. -- P.~291--294.

\item
Hencl S., Koskela P. {\it Regularity of the inverse of a planar Sobolev homeomorphism.}//
Arch. Rational Mech. Anal. -- 2006. -- V.~180. -- N.~1. -- P.~75--95.

\item
Hencl S., Koskela P., Maly Y. {\it Regularity of the inverse of a Sobolev homeomorphism in space.}//
Proc. Roy. Soc. Edinburgh Sect. A. -- 2006. -- V.~136A. -- N.~6. -- P.~1267--1285.

\item
Maz'ya V. {\it Sobolev spaces} -- Berlin: Springer Verlag. 1985.

\item
Whitney H. {\it On total differentiable and smooth functions.}//
Pacific J. Math. -- 1951. -- N.~1. -- P.~143--159.

\item
Vodop'yanov~S.~K., Ukhlov~A.~D. {\it Set functions and its applications in the theory of Lebesgue and Sobolev spaces.}//
Siberian Adv. Math. -- 2004. -- V.~14. -- N.~4. -- P.~1--48.

\item
Burenkov V.~I. {\it Sobolev Spaces on Domains}- Stuttgart: Teubner-Texter zur Mathematik. 1998.

\item
Kruglikov V.~I. {\it Capacities of condensers and spatial mappings quasiconformal in the mean.}//
Matem. sborn. -- 1986. -- V.~130. -- N.~2. -- P.~185--206.

\item
Rado T.,Reichelderfer P.~V. {\it Continuous Transformations in Analysis} -- Berlin: Sp\-rin\-ger Verlag. 1955.

\item
Gusman M. {\it Differentiation of integrals in $\mathbb R^n$} -- Moscow: Mir. 1978.

\item
Hajlasz P. {\it Change of variable formula under minimal assumptions.}//
Colloq. Math. -- 1993. -- V.~64. -- N.~1. -- P.~93--101.

\item
Reshetnyak Yu.~G. {\it Some geometrical properties of functions and mappings with generalized derivatives.}//
Siberian Math. J. -- 1966. -- V.~7. -- P.~886--919.

\item
Peshkichev Yu.~A. {\it Inverse mappings for homeomorphisms of the class $BL$.}//
Mathematical Notes -- 1993. -- V.~53. -- N.~5. -- P.~98--101.

\item
Federer H. {\it Geometric measure theory} -- Berlin: Sp\-rin\-ger Verlag. 1969.

\item
Hajlasz P. {\it Sobolev mappings, co-area formula and related topics.}//
Proc. on analysis and geometry. Novosibirsk -- 2000. -- P.~227--254.

\end{enumerate}
\end{document}